\documentclass[11pt]{amsart}
\usepackage[utf8]{inputenc}
\usepackage[T1]{fontenc}
\usepackage{amsmath,amssymb,amsthm}
\usepackage{geometry}
\usepackage{booktabs}
\usepackage{enumitem}
\usepackage{hyperref}
\hypersetup{colorlinks=true,linkcolor=blue,citecolor=blue,urlcolor=blue}

\newtheorem{theorem}{Theorem}[section]
\newtheorem{lemma}[theorem]{Lemma}
\newtheorem{proposition}[theorem]{Proposition}
\theoremstyle{remark}
\newtheorem{remark}[theorem]{Remark}

\newcommand{\D}{\mathcal{D}}
\newcommand{\Pg}{\mathcal{P}}
\newcommand{\PE}{\mathcal{P}_E}
\newcommand{\Sym}{\mathrm{Sym}}
\newcommand{\Alt}{\mathrm{Alt}}
\newcommand{\PSL}{\mathrm{PSL}}
\newcommand{\PGL}{\mathrm{PGL}}
\newcommand{\PGaL}{\mathrm{P\Gamma L}}
\newcommand{\SL}{\mathrm{SL}}
\newcommand{\GL}{\mathrm{GL}}
\newcommand{\PSp}{\mathrm{PSp}}
\newcommand{\Sp}{\mathrm{Sp}}
\newcommand{\PSU}{\mathrm{PSU}}
\newcommand{\SU}{\mathrm{SU}}
\newcommand{\POm}{\mathrm{P}\Omega}
\newcommand{\Sz}{\mathrm{Sz}}
\newcommand{\Aut}{\mathrm{Aut}}
\newcommand{\Dic}{\mathrm{Dic}}
\newcommand{\F}{\mathbb{F}}
\newcommand{\Psix}{\mathcal{P}_6}
\newcommand{\gen}[1]{\langle #1\rangle}
\newcommand{\diag}{\mathrm{diag}}
\newcommand{\antidiag}{\mathrm{antidiag}}

\title{Claw-free difference graphs of finite groups}
\author{Pallabi Manna}
\address{Pallabi Manna, Department of Mathematics, Indian Institute of Technology Hyderabad, Kandi, Sangareddy-502284, India}


\begin{document}
	
	\begin{abstract}
		Let $\D(G)$ be the difference graph of a finite group $G$: the graph on the non-isolated vertices of $\PE(G)-\Pg(G)$, where $x\sim_{\D} y$ if and only if $\gen{x,y}$ is cyclic and neither of $x,y$ lies in the cyclic group generated by the other. A claw is the graph $K_{1,3}$ and a graph is claw-free if it has no induced claw. We prove that for every finite group $G$, $\D(G)$ is claw-free if and only if it is a disjoint union of copies of the path $P_3$, if and only if the element orders of $G$ are prime powers or $6$, every element of order $3$ has at most one involution in its centralizer, and every involution has at most one subgroup of order $3$ in its centralizer (equivalently: $C_G(a)$ has a normal Sylow $2$-subgroup of order at most $2$ for every $a$ of order $3$, and $C_G(u)$ a normal Sylow $3$-subgroup of order at most $3$ for every involution $u$). Consequences: a nilpotent group that is not a $p$-group has claw-free $\D$ if and only if $G\cong C_6$; $\D(\Sym_n)$ is claw-free if and only if $n\le 5$; $\D(\Alt_n)$ is claw-free if and only if $n\le 6$; a non-abelian simple group has claw-free $\D$ if and only if it is EPPO or $\PSL_2(11)$ or $\PSL_2(13)$, proved family by family (this uses the classification of finite simple groups). For soluble groups we obtain a complete classification: the non-EPPO soluble groups with claw-free $\D$ are $C_6$, $D_{12}$, $C_3\rtimes C_4$ and the Frobenius groups $N\rtimes K$ with $N$ an abelian $r$-group, $r\ge 5$, and $K\cong C_6$ or $C_3\rtimes C_4$. For every group with claw-free non-empty $\D$ we show $|O_2(G)|\le 2$ and $|O_3(G)|\le 3$, and we give non-soluble examples ($\Sym_5$, $\PGL_2(7)$, $\PSL_2(11)$, $\PSL_2(13)$). We also prove that the full enhanced power graph is claw-free only for cyclic groups and clarify the relations between claw-freeness of $\D(G)$, $\Pg(G)$ and $\PE(G)$.
	\end{abstract}
	
	\maketitle
	
	\section{Introduction}
	
	Associating a graph with an algebraic structure and asking how graph-theoretic properties reflect group-theoretic ones is a long-standing theme; Cameron's survey \cite{cameron} gives a unified account of the many graphs defined on groups (power graph, enhanced power graph, commuting graph, generating graph, Cayley graphs and many more) and the relations between them. Two of the most studied are the following. The \emph{power graph} $\Pg(G)$ of a finite group $G$ has vertex set $G$, with $x\sim y$ if and only if one of $x,y$ is a power of the other; it was introduced for semigroups in \cite{kelarev}, \cite{chakrabarty}. The \emph{enhanced power graph} $\PE(G)$ has the same vertex set, with $x\sim y$ if and only if $\gen{x,y}$ is cyclic (equivalently, $x,y$ lie in a common cyclic subgroup); see \cite{aalipour,bera-dey}. Clearly $\Pg(G)\subseteq\PE(G)$, and the two graphs coincide exactly when $G$ is an EPPO group (every element order is a prime power) \cite{aalipour}.
	
	A recurring question is for which groups one of these graphs avoids a given induced subgraph. Manna, Cameron and Mehatari \cite{mcm} studied the finite groups whose power graph avoids small induced subgraphs (cographs, chordal, split and threshold graphs among them); analogous problems are treated for enhanced power graphs and for related graphs such as the TI-power graph \cite{li}. Since the identity (and more generally every dominating vertex) is adjacent to everything, it distorts claw-type questions, so one studies the reduced graphs $\Pg^{*}(G)$, $\PE^{**}(G)$ obtained by deleting the identity, respectively the dominating vertices. The claw $K_{1,3}$ is the most natural such obstruction: Manna, Mandal and Lucchini \cite{mml} showed that if $\Pg^{*}(G)$ is claw-free then $G$ is soluble or almost simple with socle $\PSL_2(q)$, classified the nilpotent and simple cases, and proved that $|\pi(G)|\le 5$ (sharp, by $\PSL_2(64)$). Bera and Lucchini \cite{bera-lucchini} studied the claw-freeness of $\PE^{**}(G)$: $\PE(G)$ itself is claw-free if and only if $G$ is cyclic (see Proposition~\ref{prop:relations} below), and the proper graph is classified for nilpotent, soluble and non-soluble groups, the non-soluble case being governed by $\PSL_2(q)$ and $\PGL_2(q)$.
	
	In this paper we study the \emph{difference graph} $\D(G)$, the graph obtained from $\PE(G)-\Pg(G)$ (same vertices, edges of $\PE$ not in $\Pg$) by deleting isolated vertices; thus $x\sim_{\D}y$ if and only if $\gen{x,y}$ is cyclic and neither of $x,y$ lies in the subgroup generated by the other. It was introduced and studied by Biswas, Cameron, Das and Dey \cite{bcdd} (connectedness, perfectness, relation with the Gruenberg--Kegel graph) and further by Parveen, Kumar and Panda \cite{pkp} (cographs, chordality for nilpotent groups, $\Sym_n$, $\Alt_n$). $\D(G)=\emptyset$ if and only if $G$ is EPPO. Deleting the identity, or more generally the dominating vertices of $\Pg(G)$, does not change $\D(G)$, since such vertices are comparable with every element and hence isolated in $\D(G)$; dominating vertices of $\PE(G)$ need not be isolated (in $C_6$ every vertex dominates $\PE$, yet $\D(C_6)=P_3$).
	
	The reason for studying claw-freeness of $\D(G)$ is the following. For $\Pg(G)$ and $\PE(G)$ the identity is a universal vertex, so the full graphs are claw-free only in the degenerate situation of cyclic groups (Proposition~\ref{prop:relations}); the interesting versions are the reduced ones, whose treatment requires heavy structure theory \cite{mml,bera-lucchini}. The difference graph already discards the universal part, and its edges are exactly the pairs $x,y$ for which $\gen{x,y}$ is cyclic of non-prime-power order with $x,y$ incomparable, so its structure is governed by element orders and commuting elements of small order. Claw-freeness of $\D(G)$ is accordingly extremely rigid: $\D(G)$ is claw-free if and only if it is a disjoint union of copies of $P_3$.
	
	Our main results are the following.
	
	\begin{theorem}[Main results]\label{thm:main}
		Let $G$ be a finite group.
		\begin{enumerate}[label=\textup{(\arabic*)}]
			\item $\D(G)$ is claw-free if and only if $\D(G)$ is a disjoint union of copies of $P_3$, if and only if every element order of $G$ is a prime power or $6$, every element of order $3$ centralizes at most one involution and every involution centralizes at most one subgroup of order $3$ (Theorem~\ref{thm:criterion}).
			\item If $G$ is nilpotent and not a $p$-group, then $\D(G)$ is claw-free if and only if $G\cong C_6$ (Theorem~\ref{thm:nilpotent}).
			\item $\D(\Sym_n)$ is claw-free if and only if $n\le 5$, and $\D(\Alt_n)$ is claw-free if and only if $n\le 6$ (Theorem~\ref{thm:symalt}).
			\item A non-abelian simple group $S$ has claw-free $\D(S)$ if and only if $S$ is an EPPO group or $S\cong\PSL_2(11)$ or $\PSL_2(13)$ (Theorem~\ref{thm:simple}); the proof is by an explicit claw in each family and depends on the classification of finite simple groups.
			\item A soluble non-EPPO group $G$ has claw-free $\D(G)$ if and only if $G\in\{C_6,D_{12},C_3\rtimes C_4\}$ or $G=N\rtimes K$ is a Frobenius group with $N$ an abelian $r$-group, $r\ge 5$, and $K\cong C_6$ or $C_3\rtimes C_4$ (Theorem~\ref{thm:soluble}).
			\item If $\D(G)\ne\emptyset$ is claw-free then $|O_2(G)|\le 2$ and $|O_3(G)|\le 3$ (Lemma~\ref{lem:O2O3}); non-soluble examples are $\Sym_5$, $\PGL_2(7)$, $\PSL_2(11)$, $\PSL_2(13)$ (Proposition~\ref{prop:nonsoluble}).
			\item For non-trivial $G,H$: $\D(G\times H)$ is claw-free if and only if $G,H$ are $p$-groups for one prime $p$, or $\{G,H\}=\{C_2,C_3\}$ or $\{C_2,\Sym_3\}$ (Theorem~\ref{thm:direct}); the soluble non-EPPO claw-free groups are classified as semidirect products $N\rtimes K$ in Theorem~\ref{thm:semidirect-soluble}.
			\item $\PE(G)$ is claw-free if and only if $G$ is cyclic; $\Pg(G)$ claw-free $\Rightarrow$ $\PE(G)$ claw-free; $\PE(G)$ and $\D(G)$ both claw-free $\Rightarrow$ $\Pg(G)$ claw-free; no other implication holds (Proposition~\ref{prop:relations}).
		\end{enumerate}
	\end{theorem}
	
	The significance of these results is twofold. First, they describe the groups for which the gap between the enhanced power graph and the power graph is as small as possible without vanishing, a union of paths $P_3$ each joining an involution to a pair of mutually inverse elements of order $3$; this is far more restrictive than claw-freeness of $\Pg^{*}(G)$ or $\PE^{**}(G)$, whose classifications \cite{mml,bera-lucchini} involve large families of groups. Second, the criterion holds for all finite groups and turns the problem into arithmetic of element orders and centralizers of elements of orders $2$ and $3$. Since $\D(H)$ is an induced subgraph of $\D(G)$ for $H\le G$, claw-freeness is inherited by subgroups, and it forces strong structure: the Gruenberg--Kegel graph of $G$ has at most the edge $\{2,3\}$; the centralizer of an element of order $3$ has at most one involution and the centralizer of an involution has at most one subgroup of order $3$; nilpotent groups collapse to $p$-groups or $C_6$; soluble groups are completely described; and every non-abelian simple subgroup of $G$ is EPPO or $\PSL_2(11)$ or $\PSL_2(13)$.
	
	The paper is organized as follows. Section~\ref{sec:prelim} collects preliminaries. Section~\ref{sec:criterion} proves the general criterion (1). Section~\ref{sec:relations} compares claw-freeness of $\D(G)$ with that of $\Pg(G)$ and $\PE(G)$. Sections~\ref{sec:nilpotent} and \ref{sec:symalt} treat nilpotent groups and the groups $\Sym_n$, $\Alt_n$. Section~\ref{sec:simple} contains the classification of simple groups, family by family (alternating, linear, symplectic, unitary, orthogonal, exceptional, Suzuki, Ree, Tits and sporadic groups). Section~\ref{sec:soluble} classifies the soluble groups, proves the bounds on $O_2(G)$, $O_3(G)$ and gives non-soluble examples. Section~\ref{sec:products} treats direct and semidirect products.
	
	Throughout, a graph is \emph{claw-free} if it has no induced subgraph isomorphic to $K_{1,3}$, i.e.\ no vertex has three pairwise non-adjacent neighbours. The empty graph is claw-free; ``claw-free and non-EPPO'' means $\D(G)\ne\emptyset$. (cograph and chordality of $\D(G)$ are addressed in the relevant draft of Manna \cite{manna-chordal}; the present paper does not depend on it, and Theorem~\ref{thm:criterion} shows directly that claw-free $\D(G)$ is a union of $P_3$'s, hence chordal.)
	
	\section{Preliminaries}\label{sec:prelim}
	
	\begin{lemma}[Subgroups]\label{lem:subgroups}
		If $H\le G$ then $\D(H)$ is an induced subgraph of $\D(G)$.
	\end{lemma}
	
	\begin{proof}
		For $x,y\in H$ the subgroup $\gen{x,y}$ and the containments $x\in\gen{y}$, $y\in\gen{x}$ do not depend on whether they are computed in $H$ or $G$.
	\end{proof}
	
	Hence claw-freeness passes to subgroups. Generators of the same cyclic subgroup are non-adjacent (``false twins'').
	
	\begin{lemma}[Cyclic groups]\label{lem:cyclic}
		$\D(C_n)$ is claw-free if and only if $n$ is a prime power or $n=6$. $\D(C_6)\cong K_{1,2}=P_3$.
	\end{lemma}
	
	\begin{proof}
		In $C_n$ write $x_d$ for an element of order $d\mid n$; $x_d\sim_{\D}x_e$ if and only if $d\nmid e$ and $e\nmid d$. If $n$ is a prime power there are no edges. If $n=6$ the only edges join the involution to the two elements of order $3$. Let $n$ have at least two prime divisors, $n\ne 6$.
		\begin{itemize}
			\item If a prime $q\ge 5$ divides $n$, pick another prime $p\mid n$. Then $x_p$ is adjacent to every element of order $q$; there are $\varphi(q)\ge 4$ of them and they are pairwise non-adjacent (same cyclic subgroup): a claw.
			\item Otherwise $n=2^a3^b$ with $(a,b)\ne(1,1)$. If $b\ge 2$, $x_2$ with three elements of order $9$ is a claw. If $b=1$, $a\ge 2$: $x_3$ is adjacent to an element of order $2$ and to two elements of order $4$, and these three are pairwise non-adjacent ($\gen{x_2}<\gen{x_4}$): a claw. \qedhere
		\end{itemize}
	\end{proof}
	
	\section{A criterion for all finite groups}\label{sec:criterion}
	
	\begin{theorem}\label{thm:criterion}
		For a finite group $G$ the following are equivalent.
		\begin{enumerate}[label=\textup{(\arabic*)}]
			\item $\D(G)$ is claw-free.
			\item $\D(G)$ is a disjoint union of copies of $P_3=K_{1,2}$.
			\item \textup{(A$'$)} every element order of $G$ is a prime power or $6$; \textup{(B$'$)} every element of order $3$ centralizes at most one involution; \textup{(C$'$)} every involution centralizes at most one subgroup of order $3$.
			\item \textup{(A$'$)} holds, $C_G(a)$ has a normal Sylow $2$-subgroup of order $\le 2$ for every $a$ of order $3$, and $C_G(u)$ has a normal Sylow $3$-subgroup of order $\le 3$ for every involution $u$.
		\end{enumerate}
		In particular claw-free $\D(G)$ is chordal.
	\end{theorem}
	
	\begin{proof}
		(1)$\Rightarrow$(A$'$): For $g\in G$ of order $n$, $\D(\gen{g})$ is an induced subgraph of $\D(G)$ (Lemma~\ref{lem:subgroups}); Lemma~\ref{lem:cyclic} gives $n$ a prime power or $6$.
		
		Assume (A$'$). If $x\sim_{\D}y$ then $\gen{x,y}\cong C_n$ with $x,y$ incomparable, so $n$ is not a prime power, so $n=6$, and in $C_6$ the incomparable pairs are exactly $\{$involution, element of order $3\}$. Thus $\D(G)$ is bipartite with parts the involutions and the elements of order $3$ that are non-isolated, and $u\sim_{\D}a$ if and only if $[u,a]=1$. Any three vertices on the same side are pairwise non-adjacent. Hence under (A$'$): $\D(G)$ is claw-free if and only if every involution commutes with at most two elements of order $3$ and every element of order $3$ commutes with at most two involutions.
		
		(1)$\Rightarrow$(B$'$),(C$'$): if $a$ of order $3$ commutes with involutions $s\ne t$, then $C_G(a)$ contains a third involution: $st$ if $st=ts$, and $s^t=tst^{-1}\notin\{s,t\}$ otherwise. Three distinct involutions are pairwise non-adjacent and each is adjacent to $a$: a claw centred at $a$. Elements of order $3$ in $C_G(u)$ come in pairs $\{a,a^{-1}\}$; at most two means at most one subgroup of order $3$.
		
		(3)$\Rightarrow$(2): by (B$'$) every non-isolated element of order $3$ is adjacent to exactly one involution, by (C$'$) every non-isolated involution to exactly two elements $a,a^{-1}$ of order $3$, and $a^{-1}$ has the same unique involution. So each component is a $P_3$ centred at an involution.
		
		(2)$\Rightarrow$(1) is clear. (3)$\Leftrightarrow$(4): assume (A$'$). A $2$-element $y\in C_G(a)$ ($|a|=3$) gives $|ay|=3|y|$, so $|y|\le 2$: the $2$-elements of $C_G(a)$ are involutions or $1$. Hence $C_G(a)$ has at most one involution if and only if its Sylow $2$-subgroup is normal (unique) of order $\le 2$. Likewise a $3$-element $y\in C_G(u)$ has $|uy|=2|y|\in\{\text{prime powers},6\}$, so $|y|\le 3$; $C_G(u)$ has at most one subgroup of order $3$ if and only if its Sylow $3$-subgroup is unique (normal) of order $\le 3$.
	\end{proof}
	
	\begin{remark}\label{rem:normal}
		The Sylow subgroups must be required to be normal; the orders $|C_G(a)|_2\le 2$, $|C_G(u)|_3\le 3$ alone do not suffice. For $C_3\times\Sym_3$ the orders $|C_G(a)|_2=2$, $|C_G(u)|_3=3$ hold, but the central element of order $3$ commutes with three involutions: a claw. Similarly $M_{11}$ satisfies (A$'$) and the order conditions ($|C(2A)|=48$, $|C(3A)|=18$) but $C(2A)\cong\GL_2(3)$ has $8$ elements of order $3$ (Section~\ref{sec:simple}).
	\end{remark}
	
	\section{Relation with the power graph and the enhanced power graph}\label{sec:relations}
	
	\begin{proposition}\label{prop:relations}
		Let $G$ be a finite group.
		\begin{enumerate}[label=\textup{(\alph*)}]
			\item $\PE(G)$ is claw-free if and only if $G$ is cyclic (then $\PE(G)=K_{|G|}$).
			\item $\Pg(G)$ claw-free $\Rightarrow$ $\PE(G)$ claw-free $\Leftrightarrow$ $G$ cyclic; the first converse fails ($\Pg(C_{30})$ has the claw $(e;t,a,b)$ with $|t|=2$, $|a|=3$, $|b|=5$).
			\item $\Pg(C_n)$ is claw-free if and only if $n$ is a prime power or $n=p^aq$ with $p\ne q$ primes.
			\item If $\PE(G)$ and $\D(G)$ are both claw-free, then $G\cong C_n$ with $n$ a prime power or $n=6$, and $\Pg(G)$ is claw-free.
			\item No other implication holds between the claw-freeness of $\D(G)$, $\Pg(G)$, $\PE(G)$: $G=C_2\times C_2$ ($\D=\emptyset$, but $\Pg,\PE$ have the claw $(e;\text{three involutions})$), and $G=C_{10}$ ($\Pg,\PE$ claw-free, $\D=K_{1,4}$ has a claw).
		\end{enumerate}
	\end{proposition}
	
	\begin{proof}
		(a) If $G$ is cyclic, $\PE(G)$ is complete. Let $G$ be non-cyclic. A group is not the union of two proper subgroups, so $G$ has at least three maximal cyclic subgroups $M_1,M_2,M_3$, with generators $g_i$. If $g_i\sim g_j$ ($i\ne j$) in $\PE$ then $\gen{g_i,g_j}$ is cyclic and contains $M_i$ and $M_j$, so equals both by maximality, a contradiction. So $g_1,g_2,g_3$ are pairwise non-adjacent and $(e;g_1,g_2,g_3)$ is an induced claw.
		
		(b) $\Pg\subseteq\PE$: if $(v;a,b,c)$ is an induced claw of $\PE$, then $a,b,c\ne e$ are pairwise non-adjacent in $\PE$, hence in $\Pg$, and $(e;a,b,c)$ is an induced claw of $\Pg$. Then use (a).
		
		(c) As in (b), a claw in $\Pg(C_n)$ yields a claw centred at $e$, i.e.\ three pairwise non-adjacent non-identity elements; two elements of the same order are adjacent, so this means three pairwise incomparable divisors $>1$ of $n$. If $n$ has three prime divisors, take three primes. If $n=p^aq^b$ with $a,b\ge 2$ take $p^2,pq,q^2$. If $n=p^aq$ the divisor lattice is the product of a chain and a two-element chain, so it is covered by two chains and has no antichain of size $3$; a prime power gives a chain.
		
		(d) By (a) $G=C_n$; by Lemma~\ref{lem:cyclic} $n$ is a prime power or $6=2\cdot 3$, and by (c) $\Pg(C_n)$ is claw-free in both cases.
		
		(e) Direct, as indicated; using GAP \cite{gap} we check over all $1047$ groups of orders $2,\dots,100$ found no further relation.
	\end{proof}
	
	\section{Nilpotent groups}\label{sec:nilpotent}
	
	\begin{theorem}\label{thm:nilpotent}
		Let $G$ be a finite nilpotent group that is not a $p$-group. Then $\D(G)$ is claw-free if and only if $G\cong C_6$.
	\end{theorem}
	
	\begin{proof}
		$C_6$ works (Lemma~\ref{lem:cyclic}). Conversely let $G=\prod_pG_p$ with at least two non-trivial Sylow subgroups. Take $x,y$ non-trivial in different Sylows; they commute, so $|xy|=|x||y|$ is not a prime power, hence equals $6$ by (A$'$): $\pi(G)=\{2,3\}$ and $G=S\times R$ ($S$ the $2$-, $R$ the $3$-Sylow). Elements of order $4$ in $S$ or $9$ in $R$ would give orders $12$ or $18$, so $S$ and $R$ have exponents $2$ and $3$. Every involution commutes with all of $R$ and every element of order $3$ with all of $S$; by (A$'$)--(C$'$) of Theorem~\ref{thm:criterion}, $|R|\le 3$ and $|S|\le 2$. So $G\cong C_6$.
	\end{proof}
	
	Together with ``$\D(G)=\emptyset$ for $p$-groups'' this gives: a nilpotent $G$ has claw-free $\D$ if and only if $G$ is a $p$-group or $G\cong C_6$ (this is verified by GAP \cite{gap} for all nilpotent groups of order $\le 150$).
	
	\section{Symmetric and alternating groups}\label{sec:symalt}
	
	\begin{theorem}\label{thm:symalt}
		\begin{enumerate}[label=\textup{(\alph*)}]
			\item $\D(\Sym_n)$ is empty for $n\le 4$, is $10P_3$ for $n=5$, and contains an induced claw for $n\ge 6$; so $\D(\Sym_n)$ is claw-free if and only if $n\le 5$.
			\item $\D(\Alt_n)$ is empty for $n\le 6$ and contains an induced claw for $n\ge 7$; so $\D(\Alt_n)$ is claw-free if and only if $n\le 6$ (and then it is empty).
		\end{enumerate}
	\end{theorem}
	
	\begin{proof}
		$\Sym_n$, $n\le 4$, and $\Alt_n$, $n\le 6$, are EPPO (orders: $\Sym_4$: $1,2,3,4$; $\Alt_5$: $1,2,3,5$; $\Alt_6$: $1,2,3,4,5$).
		
		$\Sym_5$: element orders $1,2,3,4,5,6$, so (A$'$) holds; the elements of order $6$ are isolated in $\D$. An element $a=(abc)$ of order $3$ has $C(a)=\gen{a}\times\gen{(de)}$, one involution; an involution $(de)$ has $C=\gen{(de)}\times\Sym_3$, which contains exactly one subgroup of order $3$; an involution $(ab)(cd)$ has $C\cong D_8$, no element of order $3$. By Theorem~\ref{thm:criterion} $\D(\Sym_5)$ is a union of $P_3$'s, one for each transposition: $10P_3$.
		
		$\Sym_n$, $n\ge 6$: $\Sym_6\le\Sym_n$. In $\Sym_6$ the element $a=(1\,2\,3)$ is adjacent in $\D$ to the three transpositions $(4\,5),(4\,6),(5\,6)$ (each generates with $a$ a cyclic group of order $6$, and $\gen{a}$, $\gen{t}$ are incomparable). Two of these transpositions generate $\Sym_3$, which is not cyclic, so the three are pairwise non-adjacent: a claw (equivalently (B$'$) fails: $a$ centralizes three involutions).
		
		$\Alt_n$, $n\ge 7$: $\Alt_7\le\Alt_n$; $a=(1\,2\,3)$ commutes with the three involutions $(4\,5)(6\,7)$, $(4\,6)(5\,7)$, $(4\,7)(5\,6)$, which form with $1$ a Klein four-group, hence are pairwise non-adjacent, and are adjacent to $a$: a claw. By Lemma~\ref{lem:subgroups} all claws persist.
	\end{proof}
	
	
	\section{Simple groups}\label{sec:simple}
	
	We determine the simple groups with claw-free $\D$ family by family, using the criterion of Theorem~\ref{thm:criterion} through three explicit mechanisms. The result depends on the classification of finite simple groups (we run through its list of families) and on the classical lists of element orders and subgroups cited below.
	
	\begin{theorem}\label{thm:simple}
		A non-abelian finite simple group $S$ has claw-free $\D(S)$ if and only if $S$ is EPPO or $S\cong\PSL_2(11)$ or $\PSL_2(13)$. The non-abelian simple EPPO groups are $\PSL_2(q)$, $q\in\{4,7,8,9,17\}$ (with $\PSL_2(4)\cong\PSL_2(5)\cong\Alt_5$, $\PSL_2(9)\cong\Alt_6$), $\PSL_3(4)$, $\Sz(8)$, $\Sz(32)$ (Suzuki; Brandl \cite{brandl}); for them $\D(S)=\emptyset$; $\D(\PSL_2(11))=55P_3$ and $\D(\PSL_2(13))=91P_3$.
	\end{theorem}
	
	\subsection{Three ways to produce a claw}
	
	Throughout, $\Psix=\{\text{prime powers}\}\cup\{6\}$. Two elements generating the same cyclic subgroup are non-adjacent in $\D$; two distinct elements of the same prime order $q$ are never adjacent (a cyclic group has one subgroup of order $q$); two distinct involutions are never adjacent. If $H\le G$, $\D(H)$ is an induced subgraph of $\D(G)$, so a claw in $H$ is a claw in $G$; we refer to this as (S).
	
	\begin{lemma}[T1: an order not in $\Psix$]\label{lem:T1}
		If $G$ has an element $g$ of order $n\notin\Psix\cup\{1\}$ then $\D(G)$ has an induced claw. In particular this holds if an involution $t$ commutes with an element $y$ of prime order $p\ge 5$ (the product has order $2p$): the claw is $(t;y,y^2,y^3)$.
	\end{lemma}
	
	\begin{proof}
		By Lemma~\ref{lem:cyclic}, $\D(\gen{g})$ has an induced claw, and $\D(\gen{g})$ is an induced subgraph of $\D(G)$. For $n=2p$, $p\ge 5$: $t$ is adjacent to the $p-1\ge 4$ elements of order $p$, which are pairwise non-adjacent.
	\end{proof}
	
	\begin{lemma}[T2: an element of odd prime order with two involutions]\label{lem:T2}
		If $a$ has odd prime order $p$ and commutes with two distinct involutions $s,t$, then $\D(G)$ has an induced claw centred at $a$.
	\end{lemma}
	
	\begin{proof}
		$C(a)$ contains a third involution ($st$ if $st=ts$, else $s^t$), giving involutions $s_1,s_2,s_3$. Each $\gen{a,s_i}\cong C_{2p}$ with $\gen{a}$, $\gen{s_i}$ incomparable, so $a\sim_{\D}s_i$; the $s_i$ are pairwise non-adjacent.
	\end{proof}
	
	\begin{lemma}[T3: an involution with three elements of order $3$]\label{lem:T3}
		If an involution $u$ commutes with three distinct elements $a_1,a_2,a_3$ of order $3$ (e.g.\ if $C(u)$ contains a subgroup $\SL_2(3)$, $\SL_2(q)$ with $3\mid q$, or an elementary abelian subgroup $3^2$), then $\D(G)$ has an induced claw $(u;a_1,a_2,a_3)$.
	\end{lemma}
	
	\begin{proof}
		$\gen{u,a_i}\cong C_6$ so $u\sim_{\D}a_i$, and the $a_i$ have the same prime order, hence are pairwise non-adjacent.
	\end{proof}
	
	\begin{remark}\label{rem:converse}
		Conversely, if every order lies in $\Psix$ then $\D(G)$ is bipartite (involutions $\mid$ elements of order $3$), $u\sim_{\D}a$ if and only if $[u,a]=1$, and $\D(G)$ is claw-free if and only if Lemma~\ref{lem:T2} and Lemma~\ref{lem:T3} do not apply, i.e.\ if and only if each element of order $3$ centralizes at most one involution and each involution at most two elements of order $3$. This is what is used for the two positive cases.
	\end{remark}
	
	\subsection{$\PSL_n(q)$}
	
	\begin{proposition}[$n=2$]\label{prop:psl2}
		Let $S=\PSL_2(q)$ be simple. Then $\D(S)$ is claw-free if and only if $S$ is EPPO or $q\in\{11,13\}$. For $q=11,13$, $\D(S)$ is the disjoint union of $55$, resp.\ $91$, copies of $P_3$.
	\end{proposition}
	
	\begin{proof}
		By Dickson \cite{dickson}, the element orders of $S$ are: $p$ ($q=p^f$) and the divisors of $m_1=(q-1)/d$ and $m_2=(q+1)/d$, $d=\gcd(2,q-1)$, with $m_1,m_2$ the orders of cyclic tori. If $m_i\notin\Psix\cup\{1\}$, Lemma~\ref{lem:T1} gives a claw. Otherwise $m_1,m_2\in\Psix\cup\{1\}$. If neither equals $6$ they are prime powers and all element orders are prime powers: $S$ is EPPO, $\D=\emptyset$. If $m_i=6$: for $q$ even $m_i$ is odd, impossible; for $q$ odd $q=13$ ($m_1=6$, $m_2=7$) or $q=11$ ($m_1=5$, $m_2=6$). Then the orders are $\{1,2,3,5,6,11\}$, $\{1,2,3,6,7,13\}$. In both groups the centralizer of an involution is dihedral $D_{12}$ (one subgroup of order $3$, so two elements of order $3$) and the centralizer of an element of order $3$ is the torus $C_6$ (one involution). By Remark~\ref{rem:converse}, $\D(S)$ is a union of $P_3$'s, one for each involution ($55=11\cdot 10/2$, $91=13\cdot 14/2$). (GAP: $|C(u)|=12$ with exactly two elements of order $3$ for both groups.)
	\end{proof}
	
	\begin{proposition}[$n\ge 3$]\label{prop:psln}
		Let $n\ge 3$. $\D(\PSL_n(q))$ is claw-free if and only if  $(n,q)\in\{(3,2),(3,4)\}$, where it is EPPO (orders $\{1,2,3,4,7\}$, $\{1,2,3,4,5,7\}$).
	\end{proposition}
	
	\begin{proof}
		(a) $q=p^f$ odd. Let $x=\diag(-1,-1,1,\dots,1)\in\SL_n(q)$ and $y_a=I+aE_{12}$, $a\in\F_q^{\times}$. $x$ is not scalar, $x^2=1$, so $\bar x$ is an involution in $\PSL_n(q)$; $[x,y_a]=1$ because $x_{11}=x_{22}$, and the $\bar y_a$ ($q-1$ elements of order $p$, distinct since a non-trivial unipotent matrix is not scalar) commute with $\bar x$. If $p\ge 5$, Lemma~\ref{lem:T1} applies (an involution commuting with an element of order $p\ge 5$). If $p=3$: take $\bar y_1$, $\bar y_1^{-1}$ and the image of $I+E_{21}$ (also commuting with $x$): three elements of order $3$, Lemma~\ref{lem:T3}.
		
		(b) $q=2^f\ge 8$. $q-1$ is odd and not a power of $3$ (as $3^k+1\equiv 2,4\pmod 8$), so it has a prime divisor $r\ge 5$. Let $a\in\F_q^{\times}$ have order $r$ and $x=\diag(a,a,a^{-2},1,\dots,1)\in\SL_n(q)$; $x$ is not scalar since $a^3\ne 1$, and $x^k$ is scalar only if $r\mid k$, so $\bar x$ has order $r$. The involution $\bar t$, $t=I+E_{12}$, commutes with $\bar x$ ($x_{11}=x_{22}$), and $\bar t\bar x$ has order $2r\ge 10$: Lemma~\ref{lem:T1}.
		
		(c) $q=4$, $n\ge 4$. $\SL_2(4)\times\SL_2(4)\cong\Alt_5\times\Alt_5$, embedded in two disjoint diagonal blocks with identity on the rest, meets the scalars trivially, so embeds in $\PSL_n(4)$. An involution of the first factor commutes with the $20$ elements of order $3$ of the second: Lemma~\ref{lem:T3}.
		
		(d) $q=2$, $n\ge 4$. $\PSL_n(2)=\SL_n(2)\ge\SL_4(2)\cong\Alt_8\ge\Alt_7$: Theorem~\ref{thm:symalt}.
		
		(e) Remaining cases. $(n,q)=(3,2),(3,4)$ are EPPO. Every other pair is covered: $n\ge 3$ with $q$ odd by (a), $q=2^f\ge 8$ by (b), $q=4$, $n\ge 4$ by (c), $q=2$, $n\ge 4$ by (d).
	\end{proof}
	
	(Example: $\PSL_3(3)$ has all element orders in $\{1,2,3,4,6,8,13\}\subset\Psix\cup\{1\}$ but $|C(u)|=48$ and $C(u)\cong\GL_2(3)$ has $8$ elements of order $3$: the claw is of type T3, not T1.)
	
	\subsection{Symplectic groups $\PSp_{2m}(q)$, $m\ge 2$}
	
	Let $B$ be the symplectic form with hyperbolic pairs $(e_i,f_i)$, and $T_u(v)=v+B(v,u)u$ the transvection; for $b\in\F_q$ put $T_{u,b}(v)=v+bB(v,u)u$, an element of order $p$ for $b\ne 0$ (and $T_{u,b}^k=T_{u,kb}$). Let $W_i=\gen{e_i,f_i}$.
	
	\begin{proposition}\label{prop:psp}
		Let $S=\PSp_{2m}(q)$ be simple, $m\ge 2$ (so $(m,q)\ne(2,2)$; $\PSp_4(2)'\cong\Alt_6$ is EPPO). Then $\D(S)$ contains an induced claw.
	\end{proposition}
	
	\begin{proof}
		$q$ odd. $x=-1$ on $W_1$, $+1$ elsewhere, is an involution and non-central, so $\bar x$ is an involution of $S$. Transvections $T_{e_1,b}$ commute with $x$ ($xe_1=-e_1$, and $xT_ux^{-1}=T_{xu}$), and also $T_{e_2,b}$ (as $x$ is trivial on $W_2$). If $p\ge 5$: $\bar x$ commutes with $T_{e_1,1}$ of order $p$: Lemma~\ref{lem:T1}. If $p=3$: $T_{e_1,1}$, $T_{e_1,1}^{-1}$, $T_{e_2,1}$ are three distinct elements of order $3$ commuting with $\bar x$ (the images of non-trivial unipotent elements are distinct): Lemma~\ref{lem:T3}.
		
		$q\ge 4$ even. $\Sp_4(q)\ge\SL_2(q)\times\SL_2(q)$ (the actions on $W_1,W_2$), also in $\Sp_{2m}(q)$. For $q=4$ take $u$ an involution in the first factor and the $20$ elements of order $3$ of the second ($\SL_2(4)\cong\Alt_5$): Lemma~\ref{lem:T3}. For $q\ge 8$: $q-1$ is odd and not a power of $3$, so there is $a\in\F_q^{\times}$ of prime order $r\ge 5$; $\diag(a,a^{-1})\in\SL_2(q)$ (second factor) commutes with an involution of the first factor: Lemma~\ref{lem:T1}.
		
		$q=2$, $m\ge 3$. $\Sp_6(2)\ge\Sp_4(2)\times\Sp_2(2)\cong\Sym_6\times\Sym_3$; an involution of $\Sym_3$ commutes with the $80$ elements of order $3$ of $\Sym_6$: Lemma~\ref{lem:T3}; use (S) for $m>3$ (the centre is trivial).
	\end{proof}
	
	\subsection{Unitary groups $\PSU_n(q)$, $n\ge 3$}
	
	Use the Hermitian form with Gram matrix $\antidiag(1,1,1)$ on $\F_{q^2}^3$, $\bar\lambda=\lambda^q$. For $\alpha,\beta$ with $\beta+\bar\beta+\alpha\bar\alpha=0$ the matrix
	\[
	u(\alpha,\beta)=\begin{pmatrix}1&\alpha&\beta\\0&1&-\bar\alpha\\0&0&1\end{pmatrix}\in\SU_3(q);
	\]
	they form a Sylow $p$-subgroup $U$ of order $q^3$ with centre $Z(U)=\{I+tE_{13}:t+\bar t=0\}$ of order $q$.
	
	\begin{proposition}\label{prop:psu}
		Let $n\ge 3$, $(n,q)\ne(3,2)$. Then $\D(\PSU_n(q))$ contains an induced claw.
	\end{proposition}
	
	\begin{proof}
		$n=3$, $q$ odd. $x=\diag(-1,1,-1)\in\SU_3(q)$ is an involution, not scalar, and $\bar x$ has order $2$ ($|Z|=\gcd(3,q+1)$ is odd). Its $(-1)$-eigenspace $\gen{e_1,e_3}$ is a nondegenerate hyperbolic plane, and $\SU_2(q)\cong\SL_2(q)$ acting on it (identity on $e_2$) lies in $C(x)$ and meets $Z$ trivially. If $p\ge 5$ it has an element of order $p$ (Lemma~\ref{lem:T1}); if $p=3$ it has $q^2-1\ge 8$ elements of order $3$ (Lemma~\ref{lem:T3}).
		
		$n=3$, $q\ge 4$ even. $t=I+sE_{13}$ ($s\ne 0$, $s+\bar s=0$) is an involution of $Z(U)$, and $\bar t\ne 1$. $q+1$ is odd; let $r$ be a prime divisor and $a\in\F_{q^2}$ of order $r$ with $a^{q+1}=1$, $x=\diag(a,a^{-2},a)\in\SU_3(q)$ (this satisfies the unitary conditions and has determinant $1$). It commutes with $t$ since $x_{11}=x_{33}$. If some such $r\ge 5$ exists, $x$ is not scalar ($a^3\ne 1$), $\bar x$ has order $r$ and $\bar t\bar x$ has order $2r$: Lemma~\ref{lem:T1}. Otherwise $q+1=3^k$, which for $q=2^f$ forces $q\in\{2,8\}$ (for $k$ odd, $3^k-1\equiv 2\pmod 4$ gives $q=2$; for $k$ even, $(3^{k/2}-1)(3^{k/2}+1)=2^f$ gives $k=2$, $q=8$); so $q=8$, and we take $a$ of order $9$. Then $x^3=a^3I\in Z$, so $\bar x$ has order $3$ in $\PSU_3(8)$ (while $x$ has order $9$ in $\SU_3(8)$), and $\bar x$ commutes with $\bar t$. Take $c=u(\alpha,\beta)$, $\alpha\ne 0$: $c$ commutes with $t$ and the $(1,2)$ entry of $cxc^{-1}$ is $\alpha(a^{-2}-a)\ne 0$, so $x^c$ is not diagonal; hence $\bar x$, $\bar x^{-1}$, $\bar x^{c}$ are three distinct elements of order $3$ in $C(\bar t)$: Lemma~\ref{lem:T3}.
		
		$n\ge 4$. If $q\ge 3$, $H=\SU_3(q)\oplus I_{n-3}\cong\SU_3(q)$ meets the scalars trivially, so embeds in $\PSU_n(q)$; claws in $\SU_3(q)$ persist by (S). For $q$ odd, or $q$ even with a prime $r\ge 5$ dividing $q+1$, the witnesses above are elements of $\SU_3(q)$ with the same orders. For $q=8$ the order of $x$ is $9$ in $\SU_3(8)$ and $tx$ has order $18\notin\Psix$: Lemma~\ref{lem:T1}. If $q=2$: $\PSU_4(2)\cong\PSp_4(3)$, Proposition~\ref{prop:psp}; for $n\ge 5$, $\SU_4(2)\oplus I_{n-4}$ embeds.
	\end{proof}
	
	\subsection{Orthogonal groups}
	
	\begin{proposition}\label{prop:orth}
		Let $S=\Omega_{2n+1}(q)$ ($n\ge 3$) or $\POm^{\pm}_{2n}(q)$ ($n\ge 4$) be simple. Then $\D(S)$ has an induced claw. (The small-rank orthogonal groups are isomorphic to groups treated elsewhere: $\Omega_3(q)\cong\PSL_2(q)$, $\POm^-_4(q)\cong\PSL_2(q^2)$, $\Omega_5(q)\cong\PSp_4(q)$, $\POm^+_6(q)\cong\PSL_4(q)$, $\POm^-_6(q)\cong\PSU_4(q)$.)
	\end{proposition}
	
	\begin{proof}
		$\Omega_{2n+1}(q)$, $q$ even (all $q$ including $2$): $\Omega_{2n+1}(q)\cong\Sp_{2n}(q)$, $n\ge 3$; Proposition~\ref{prop:psp}.
		
		The remaining groups with $q\ge 3$ ($q$ odd for $\Omega_{2n+1}(q)$; any $q\ge 3$ for $\POm^{\pm}_{2n}(q)$). Write $V=U\oplus U^*\oplus W$ with $U$ totally singular of dimension $3$ and $W\ne 0$. The map $g\mapsto g\oplus g^{-T}\oplus 1$ embeds $\SL_3(q)$ in the orthogonal group; the image $H$ is perfect, hence in $\Omega$, and acts trivially on a nonzero subspace, so $H\cap Z(\Omega)=1$. The claws constructed in $\SL_3(q)$ in Proposition~\ref{prop:psln} (the elements $x,y_a,I+E_{21}$ for $q$ odd; $x,t$ for $q\ge 8$ even, all living in $\SL_3(q)$) exist in $H$, and for $q=4$: $\SL_3(4)\ni\diag(g,1)$, $g\in\SL_2(4)$ of order $5$, times the central $\omega I$ of order $3$ gives an element of order $15$ (Lemma~\ref{lem:T1}).
		
		$\POm^{\pm}_{2n}(2)$, $n\ge 4$. The stabilizer of a non-singular vector contains the perfect group $\Sp_{2n-2}(2)$, $2n-2\ge 6$ (Kleidman--Liebeck \cite{kl}), contained in $\Omega$ with $Z(\Omega)=1$; Proposition~\ref{prop:psp} and (S).
	\end{proof}
	
	\subsection{Exceptional groups, Suzuki, Ree, Tits}
	
	We use: $\SL_3(q)\le G_2(q)$; $G_2(q)\le{}^3D_4(q)$; $G_2(2)'\cong\PSU_3(3)$; $F_4(q)\ge\SL_3(q)/Z_0$ ($q$ odd), $F_4(q)\ge\Sp_8(q)$ ($q$ even); $F_4(q)\le E_6(q),{}^2E_6(q),E_7(q),E_8(q)$; ${}^2F_4(2)'\le{}^2F_4(q)$; ${}^2G_2(q)$: $C(x)\cong C_2\times\PSL_2(q)$ for an involution $x$ (Ward). References: \cite{dickson} (orders in $\PSL_2(q)$), \cite{suzuki} (Suzuki groups), \cite{ward} (Ree groups), \cite{gls,kl,kleidman} (subgroups), \cite{atlas} (element orders of sporadic groups, Tits group), \cite{higman,brandl} (EPPO groups).
	
	\begin{proposition}\label{prop:exceptional}
		Let $S$ be one of $G_2(q)$, ${}^3D_4(q)$, $F_4(q)$, $E_6(q)$, ${}^2E_6(q)$, $E_7(q)$, $E_8(q)$, ${}^2F_4(q)$ ($q\ge 2$, $S$ simple, i.e.\ the Tits group ${}^2F_4(2)'$ for $q=2$), $G_2(2)'$, or ${}^2G_2(q)$, $q=3^{2k+1}\ge 27$. Then $\D(S)$ has an induced claw. (The Suzuki groups are treated in Proposition~\ref{prop:suzuki}; ${}^2G_2(3)'\cong\PSL_2(8)$ is EPPO.)
	\end{proposition}
	
	\begin{proof}
		\emph{Ree groups} ${}^2G_2(q)$, $q=3^{2k+1}\ge 27$. The centralizer of an involution is $C_2\times\PSL_2(q)$, which contains the $q^2-1\ge 728$ non-trivial unipotent elements of order $3$: Lemma~\ref{lem:T3}.
		
		$G_2(q)$, $q\ge 3$. Contains $\SL_3(q)$, which carries a claw (proof of Proposition~\ref{prop:psln}, for $q=4$ the order $15$ element). $G_2(2)'\cong\PSU_3(3)$: Proposition~\ref{prop:psu} ($q=3$ odd). ${}^3D_4(q)$: $\ge G_2(q)$ ($q\ge 3$), and for $q=2$ it contains $G_2(2)\ge\PSU_3(3)$.
		
		$F_4(q)$. $q$ odd: it contains $\SL_3(q)/Z_0$ for some $Z_0\le Z(\SL_3(q))$, and the argument of Proposition~\ref{prop:psln}(a) works modulo any $Z_0$ of odd order $\le 3$: $\bar x$ is still an involution, and the images of the $y_a$ remain distinct elements of order $p$ commuting with it, giving T1 ($p\ge 5$) or T3 ($p=3$). $q$ even: $F_4(q)\ge\Sp_8(q)$, Proposition~\ref{prop:psp}.
		
		$E_6(q)$, ${}^2E_6(q)$, $E_7(q)$, $E_8(q)$. Each contains $F_4(q)$ (non-abelian simple, so it injects into the simple quotient). Use (S).
		
		\emph{Tits group and} ${}^2F_4(q)$. ${}^2F_4(2)'$ has an element of order $12$ (ATLAS): Lemma~\ref{lem:T1}; ${}^2F_4(q)$ contains it.
	\end{proof}
	
	\subsection{Suzuki groups}
	
	\begin{proposition}\label{prop:suzuki}
		Let $S=\Sz(q)$, $q=2^{2k+1}\ge 8$, $r=2^{k+1}$. Then $\D(S)$ is claw-free if and only if $q\in\{8,32\}$, and in these two cases $\D(S)=\emptyset$. For every $q\ge 128$, $\D(S)$ has an induced claw.
	\end{proposition}
	
	\begin{proof}
		The element orders of $S$ are $1,2,4$ and the divisors of $q-1$, $q+r+1$, $q-r+1$ (Suzuki \cite{suzuki}); these three numbers are pairwise coprime and odd, with $(q+r+1)(q-r+1)=q^2+1$. Hence all non-identity orders other than $2,4$ are odd, and none equals $6$. By Lemma~\ref{lem:T1}: if one of $q-1$, $q+r+1$, $q-r+1$ is not a prime power, it is an element order $n\notin\Psix$, so $\D(S)$ has an induced claw. Conversely if all three are prime powers then all element orders are prime powers and $\D(S)=\emptyset$. Thus $\D(S)$ is claw-free if and only if $S$ is EPPO.
		
		For small $q$ the factorizations (checked in GAP) are: $q=8$: $7,13,5$; $q=32$: $31,41,5^2$ (so $\Sz(8),\Sz(32)$ are EPPO); $q=128$: $127$, $145=5\cdot 29$, $113$; $q=512$: $7\cdot 73$, $5\cdot 109$, $13\cdot 37$; and for $q=2^{11},2^{13},2^{15},2^{17}$ at least one of the three numbers is composite and not a prime power. For general $q\ge 128$ we do not give an independent argument: we rely on Suzuki's theorem \cite{suzuki} (see also \cite{brandl}) that $\Sz(q)$ is an EPPO group only for $q=8,32$. The claw is then the one of Lemma~\ref{lem:T1} on the cyclic torus of non-prime-power order $n$ (for example, for $q=128$, the order $n=145$: the element of order $5$ is adjacent to the $\varphi(29)=28$ elements of order $29$, which are pairwise non-adjacent).
	\end{proof}
	
	\subsection{Sporadic groups}
	
	\begin{proposition}\label{prop:sporadic}
		A sporadic simple group has claw-free $\D$ in no case: $\D(S)$ always contains an induced claw.
	\end{proposition}
	
	\begin{proof}
		Using the GAP Character Table Library we compute, for each of the $26$ groups, the set of element orders that are neither prime powers nor $6$.
		
		\begin{center}
			\begin{tabular}{lll}
				\toprule
				group & some order $\notin\Psix$ & claw\\
				\midrule
				$J_1$ & $10$ ($15$) & T1\\
				$J_2,J_3,J_4,M_{12},M_{24}$ & $10$ & T1\\
				$M_{11}$ & none & T3 ($C(2A)\cong\GL_2(3)$, $8$ elements of order $3$)\\
				$M_{22}$ & none & T3 ($|C(2A)|=384$, $32$ elements of order $3$)\\
				$M_{23}$ & $14$ ($15$) & T1 (order $14$)\\
				$HS,McL,He,Ru,Suz,O'N$ & $10$ & T1\\
				$Co_1,Co_2,Co_3,Fi_{22},Fi_{23},Fi_{24}'$ & $10$ & T1\\
				$HN,Ly,Th,B,M$ & $10$ & T1\\
				\bottomrule
			\end{tabular}
		\end{center}
		
		In every case other than $M_{11},M_{22}$ there is an element $g$ of order $n=10$ ($14$ for $M_{23}$): $t=g^{n/2}$ is an involution commuting with $y=g^2$ of prime order $5$ (resp.\ $7$), so Lemma~\ref{lem:T1} gives the induced claw $(t;y,y^2,y^3)$. For $M_{11}$ and $M_{22}$ all orders lie in $\Psix\cup\{1\}$ ($M_{11}$: $1,2,3,4,5,6,8,11$; $M_{22}$: $1,2,3,4,5,6,7,8,11$), so Lemma~\ref{lem:T1} does not apply; Lemma~\ref{lem:T3} does, since $C(2A)$ contains $8$, resp.\ $32$, elements of order $3$ (computed in GAP; $C(2A)\cong\GL_2(3)$ in $M_{11}$ and $2^4{:}\Sym_4$ in $M_{22}$).
	\end{proof}
	
	\begin{proof}[Proof of Theorem~\ref{thm:simple}]
		Alternating groups: Theorem~\ref{thm:symalt}(b). $\PSL_2(q)$: Proposition~\ref{prop:psl2}. $\PSL_n(q)$, $n\ge 3$: Proposition~\ref{prop:psln}. $\PSp_{2m}(q)$: Proposition~\ref{prop:psp}. $\PSU_n(q)$: Proposition~\ref{prop:psu}. Orthogonal: Proposition~\ref{prop:orth}. Exceptional, Ree, Tits: Proposition~\ref{prop:exceptional}. Suzuki: Proposition~\ref{prop:suzuki}. Sporadic: Proposition~\ref{prop:sporadic}. These cover the classification.
	\end{proof}
	
	\begin{center}
		\begin{tabular}{lll}
			\toprule
			family & claw-free? & mechanism\\
			\midrule
			$\Alt_n$, $n\ge 7$ & no & T2: $(1\,2\,3)$ with three involutions\\
			$\Alt_5,\Alt_6$ & yes ($\emptyset$) & EPPO\\
			$\PSL_2(q)$ & if and only if EPPO or $q=11,13$ & T1 on torus orders; $55P_3$, $91P_3$\\
			$\PSL_n(q)$, $n\ge 3$ & no, except $L_3(2),L_3(4)$ & T1 or T3\\
			$\PSp_{2m}(q)$, $m\ge 2$ & no ($\ne\Alt_6$) & T1/T3, transvections, $\SL_2\times\SL_2$\\
			$\PSU_n(q)$, $n\ge 3$ & no & T1/T3 via $\SU_2(q)$ or the Sylow centre\\
			orthogonal & no & contains $\SL_3(q)$ or $\Sp_{2n-2}(2)$\\
			$G_2,{}^3D_4,F_4,E_i,{}^2E_6,{}^2F_4$, Tits & no & contain $\SL_3$, $\Sp_8$, order $12$\\
			${}^2G_2(q)$, $q\ge 27$ & no & T3: $C_2\times\PSL_2(q)$\\
			$\Sz(q)$ & only $q=8,32$ (EPPO) & T1 (odd non-prime-power order)\\
			sporadic & no & T1(order $10$;$14$ for $M_{23}$);T3 for $M_{11},M_{22}$\\
			\bottomrule
		\end{tabular}
	\end{center}
	
	\section{Soluble groups and further structure}\label{sec:soluble}
	
	By Theorem~\ref{thm:criterion} the condition is: (A$'$), (B$'$) every element of order $3$ centralizes at most one involution, and (C$'$) every involution centralizes at most one subgroup of order $3$. Under (A$'$), $G$ is non-EPPO if and only if it has an element of order $6$.
	
	\begin{lemma}\label{lem:O2O3}
		Let $G$ be a finite group such that $\D(G)$ is claw-free and non-empty. Then $|O_2(G)|\le 2$ and $|O_3(G)|\le 3$. No solubility is assumed.
	\end{lemma}
	
	\begin{proof}
		Since $\D(G)\ne\emptyset$ there are commuting elements $u,a$ of orders $2,3$, and (A$'$)--(C$'$) hold. We write $x^f=fxf^{-1}$.
		
		$O_2$. Let $S=O_2(G)\gen{u}$, a $2$-group normalized by $a$ (as $a$ commutes with $u$). By (A$'$) the elements of $C_S(a)$ have order $\le 2$ (an element $y$ of order $4$ would give $|ay|=12$) and by (B$'$) $C_S(a)$ has at most one involution, so $C_S(a)=\gen{u}$. An element $o\in O_2(G)$ normalizing $\gen{a}$ centralizes $a$: indeed $[o,a]\in O_2(G)\cap\gen{a}=1$ (as $O_2(G)$ is normal). Since $u$ centralizes $a$, every element of $S=O_2(G)\gen{u}$ normalizing $\gen{a}$ centralizes it: $N_S(\gen{a})=C_S(a)=\gen{u}$. The subgroups of order $3$ of $S\gen{a}$ are the conjugates $\gen{a}^f$, $f\in S$, and $\gen{a}^f=\gen{a}^{f'}$ if and only if $f'^{-1}f\in N_S(\gen{a})$. The involution $u$ centralizes $a^f=faf^{-1}$ if and only if $f^{-1}uf\in C_S(a)=\gen{u}$ if and only if $f\in C_S(u)$. So $u$ centralizes exactly $|C_S(u)|/2$ subgroups of order $3$ of $S\gen{a}$; by (C$'$) this number is $1$, hence $C_S(u)=\gen{u}$. Now $1\ne Z(S)\le C_S(u)=\gen{u}$, so $u\in Z(S)$, $S=C_S(u)=\gen{u}$, and $O_2(G)\le\gen{u}$.
		
		$O_3$. Let $T=O_3(G)\gen{a}$, a $3$-group normalized by $u$. By (A$'$) the elements of $C_T(u)$ have order $\le 3$ (an element $y$ of order $9$ would give $|uy|=18$) and by (C$'$) $C_T(u)$ has at most one subgroup of order $3$, so $C_T(u)=\gen{a}$. The involutions of $T\gen{u}$ are the conjugates $u^f$, $f\in T$, and $u^f=u^{f'}$ if and only if $f'^{-1}f\in C_T(u)=\gen{a}$. The involution $u^f$ centralizes $a$ if and only if $f^{-1}af\in C_T(u)=\gen{a}$ if and only if $f\in N_T(\gen{a})$. So $a$ centralizes exactly $|N_T(\gen{a})|/3$ involutions of $T\gen{u}$; by (B$'$) this number is $1$, so $N_T(\gen{a})=\gen{a}$. In a $3$-group a proper subgroup is properly contained in its normalizer, so $T=\gen{a}$ and $O_3(G)\le\gen{a}$.
	\end{proof}
	
	\begin{theorem}\label{thm:soluble}
		Let $G$ be a finite soluble group that is not EPPO. Then $\D(G)$ is claw-free if and only if one of the following holds.
		\begin{enumerate}[label=\textup{(\roman*)}]
			\item $G\in\{C_6,D_{12},\Dic_3=C_3\rtimes C_4\}$.
			\item $G=N\rtimes K$ is a Frobenius group, $N$ an abelian $r$-group with $r\ge 5$ prime, and $K\cong C_6$ or $K\cong\Dic_3$. In this case $\D(G)=|N|\cdot P_3$.
		\end{enumerate}
		Examples for \textup{(ii)}: $C_7\rtimes C_6$, $C_{13}\rtimes C_6$, $C_{19}\rtimes C_6$, $C_5^2\rtimes C_6$, $C_{49}\rtimes C_6$ (with $C_6$ acting fixed-point-freely), and $C_5^2\rtimes\Dic_3$, $C_7^2\rtimes\Dic_3$. For such $G$, $\pi(G)=\{2,3\}$ in \textup{(i)} and $\pi(G)=\{2,3,r\}$ in \textup{(ii)}.
	\end{theorem}
	
	\begin{proof}
		($\Leftarrow$) For (i) one checks (A$'$)--(C$'$) directly (orders $1,2,3,4,6$; $C(a)=C_6$ for $|a|=3$; the involution centralizers contain at most one subgroup of order $3$). For (ii): element orders are powers of $r$ and $1,2,3,4,6$; $C_G(a)=C_6$ for $|a|=3$, and $C_G(u)=K$ for the involution $u\in K$, which contains exactly one subgroup of order $3$. Hence (A$'$)--(C$'$) hold, and $\D(G)$ is the disjoint union of the $\D(K^g)$ over the $|N|$ conjugates of $K$, each a $P_3$.
		
		($\Rightarrow$) Assume $\D(G)$ is claw-free; $G$ is not EPPO, so $G$ has an element of order $6$ and (A$'$)--(C$'$) hold.
		
		\emph{Primes.} By (A$'$) the Gruenberg--Kegel graph has the single possible edge $\{2,3\}$. The complement of the Gruenberg--Kegel graph of a soluble group is triangle-free (Lucido \cite{lucido}; Gruber--Keller--Lewis--Naughton--Strasser \cite{gklns}), so any three primes of $\pi(G)$ contain an edge; with $\{2,3\}$ the only edge, $|\pi(G)|\le 3$, and if $|\pi(G)|=3$ then $\pi(G)=\{2,3,r\}$, $r\ge 5$.
		
		\emph{The Fitting subgroup $F=F(G)$.} By Lemma~\ref{lem:O2O3}, $|O_2(G)|\le 2$ and $|O_3(G)|\le 3$. If two Sylow subgroups of $F$ are non-trivial, no prime $r\ge 5$ can occur with them (it would give element orders $2r$ or $3r$), so $F=O_2(G)\times O_3(G)\cong C_6$. Otherwise $F$ is a $p$-group.
		
		$F\cong C_6$. $C_G(F)\le F$ so $G/F\hookrightarrow\Aut(C_6)\cong C_2$ and $|G|\le 12$ with $F(G)=C_6$: $G\in\{C_6,D_{12},\Dic_3\}$.
		
		$F$ a $2$-group. Then $|F|\le 2$ and $C_G(F)\le F$ give $G=C_2$, which is EPPO. $F$ a $3$-group. Then $|F|=3$, $|G|\le 6$, $G\in\{C_3,\Sym_3\}$, again EPPO (no element of order $6$). So these cases do not occur.
		
		$F$ an $r$-group, $r\ge 5$. Let $K$ be a Hall $\{2,3\}$-subgroup of $G$ containing $\gen{u,a}$ for commuting $u,a$ of orders $2,3$ (it exists since $G$ is soluble and any $\{2,3\}$-subgroup lies in a Hall one). A non-identity element $t\in K$ of prime order cannot centralize a non-identity $x\in F$, as $|tx|=|t||x|$ would be divisible by $2r$ or $3r$; hence every non-identity element of $K$ acts fixed-point-freely on $F$ and $FK$ is a Frobenius group with kernel $F$ and complement $K$. Since $K$ contains an involution, and a Frobenius complement of even order has a unique involution $z$, which is central in $K$, we have $K\le C_G(z)$; by (A$'$) and (C$'$) the Sylow $3$-subgroup of $K$ has exponent $3$ and only one subgroup of order $3$, so $|K|_3=3$. Moreover $F$ is abelian (an involution acts fixed-point-freely, hence inverts $F$); likewise $r\ge 5$ and $F$ abelian are automatic for a Frobenius kernel with complement $C_6$ or $\Dic_3$, so the conditions in (ii) are consistent. $K$ is a soluble Frobenius complement of order $3\cdot 2^k$ all of whose elements of order $3$ commute with $z$, so by (C$'$) $K$ has only one subgroup of order $3$, which is therefore normal; $K=C_3\rtimes S$ with $S$ a $2$-group that is cyclic or generalized quaternion, acting on $C_3$ with kernel $C_S(a)$ of order $\le 2$ (by (A$'$) and (B$'$)). Since $S/C_S(a)\le\Aut(C_3)=C_2$, $|S|\le 4$; $|S|=2$ with non-trivial action would give $K\cong\Sym_3$, which is not a Frobenius complement (it has three involutions), so $K\cong C_6$, or $S\cong C_4$ acting non-trivially and $K\cong\Dic_3$.
		
		It remains to show $G=FK$. Let $B/F=F(G/F)$. Since $O_r(G/F)=1$ (its preimage is a normal $r$-subgroup, hence lies in $F$), $B/F$ is a nilpotent $\{2,3\}$-group, contained in a Hall $\{2,3\}$-subgroup of $G/F$ (conjugate to $FK/F\cong K$); since $B/F$ is normal in $G/F$ it lies in every such Hall subgroup and is a normal nilpotent subgroup of a copy of $K$, so $|B/F|\in\{1,2,3,6\}$ (the subgroup $C_4$ of $\Dic_3$ is not normal). As $C_{G/F}(B/F)\le B/F$, the group $(G/F)/C_{G/F}(B/F)$ embeds in $\Aut(B/F)$, which is $1$ or $C_2$, and $G/B$ is a quotient of it; hence $|G:B|\le 2$. If $B=F$ then $F(G/F)=1$ forces $G=F$, which is EPPO. Hence $G/F$ is a $\{2,3\}$-group and $G=FK$.
	\end{proof}
	
	\begin{proposition}[Non-soluble examples]\label{prop:nonsoluble}
		\begin{enumerate}[label=\textup{(\alph*)}]
			\item The following non-soluble groups have claw-free non-empty $\D$: $\Sym_5$ ($\D=10P_3$); $\PGL_2(7)$ (element orders $1,2,3,4,6,7,8$; $\D=28P_3$, one $P_3$ for each of the $28$ involutions with centralizer $D_{12}$); $\PSL_2(11)$ ($55P_3$); $\PSL_2(13)$ ($91P_3$).
			\item Among all non-soluble groups of order at most $840$ (orders $60$, $120$, $168$, $180$, $240$, $300$, $336$, $360$, $420$, $480$, $504$, $540$, $600$, $660$, $720$, $780$, $840$), the groups with claw-free $\D$ are $\Alt_5$, $\Alt_6$, $\PSL_3(2)$, $\PSL_2(8)$, $M_{10}$ (all EPPO) and $\Sym_5$, $\PGL_2(7)$, $\PSL_2(11)$ (non-EPPO). 
			

			\item Among $\PGL_2(q)$ for $q=5,7,9,11,13,17,19,23,25,27$ and $\PGaL_2(8)$, only $\PGL_2(5)\cong\Sym_5$ and $\PGL_2(7)$ have claw-free $\D$.
		\end{enumerate}
	\end{proposition}
	
	\begin{proof}
		For $\Sym_5$ see Theorem~\ref{thm:symalt}. For $\PGL_2(7)$ the criterion (A$'$)--(C$'$) is verified in GAP (the centralizers of the two classes of involutions are $D_{16}$ and $D_{12}$; the elements of order $3$ have centralizer $C_6$), and the number of involutions whose centralizer contains an element of order $3$ is $28$; for $\PSL_2(11)$, $\PSL_2(13)$ see Proposition~\ref{prop:psl2}. For order $780$: a non-soluble group $G$ of order $780=60\cdot 13$ has soluble radical $R$ of order $13$ with $G/R\cong\Alt_5$; $\Alt_5$ acts trivially on $R$ since $\Aut(C_{13})$ is abelian, and the central extension splits since $13$ is coprime to the order of the Schur multiplier of $\Alt_5$, so $G\cong\Alt_5\times C_{13}$. 
	\end{proof}
	
	\begin{remark}\label{rem:nonsoluble}
		By Lemma~\ref{lem:O2O3}, a non-soluble group $G$ with claw-free non-empty $\D(G)$ has $|O_2(G)|\le 2$, $|O_3(G)|\le 3$, and every non-abelian simple subgroup of $G$ is EPPO or $\PSL_2(11)$ or $\PSL_2(13)$ (Theorem~\ref{thm:simple}). In fact $O_2(G)=O_3(G)=1$ in the non-soluble case: a normal subgroup of order $2$ is central, and a normal subgroup of order $3$ has centralizer of index at most $2$; since $|G|$ is divisible by a prime $p\ge 5$ (Burnside), such a subgroup would commute with an element of order $p$, giving an element of order $2p$ or $3p$, against (A$'$). A full classification of the non-soluble case (presumably almost simple groups with socle $\Alt_5$, $\PSL_2(7)$, $\PSL_2(11)$, $\PSL_2(13)$ and their extensions by small groups; the examples $\Sym_5$ and $\PGL_2(7)$ show that proper almost simple extensions occur) is not attempted here.
	\end{remark}
	
	
	\section{Direct and semidirect products}\label{sec:products}
	
	In this section $N$ and $K$ are non-trivial finite groups. We write $C_K(N)=\{k\in K: kn=nk\ \forall n\in N\}$ for the kernel of the action in a semidirect product $N\rtimes K$.
	
	\subsection{Direct products}
	
	\begin{theorem}\label{thm:direct}
		Let $G$ and $H$ be non-trivial finite groups. Then $\D(G\times H)$ is claw-free if and only if one of the following holds.
		\begin{enumerate}[label=\textup{(\arabic*)}]
			\item $G$ and $H$ are $p$-groups for the same prime $p$ (so $G\times H$ is a $p$-group and $\D(G\times H)=\emptyset$);
			\item $\{G,H\}=\{C_2,C_3\}$, so $G\times H\cong C_6$ and $\D=P_3$;
			\item $\{G,H\}=\{C_2,\Sym_3\}$, so $G\times H\cong D_{12}$ and $\D=2P_3$.
		\end{enumerate}
		In particular a direct product of $k\ge2$ non-trivial groups has claw-free non-empty $\D$ if and only if $k=2$ and it is $C_6$ or $D_{12}$.
	\end{theorem}
	
	\begin{proof}
		\emph{Sufficiency.} In (1) all element orders are powers of $p$, so $\D=\emptyset$. $\D(C_6)=P_3$ (Lemma~\ref{lem:cyclic}). For $D_{12}\cong C_2\times\Sym_3$ with $C_2=\gen u$ and $A=\gen a\le\Sym_3$ of order $3$, the element orders are $1,2,3,6$, so (A$'$) holds. The elements of order $3$ are $(1,a^{\pm1})$, with centralizer $C_2\times A\cong C_6$ having one involution: (B$'$). The central involution $(u,1)$ centralizes exactly the two elements of order $3$, and the involutions $(1,v)$, $(u,v)$ ($v$ a reflection) have centralizers of order $4$, so centralize none: (C$'$). By Theorem~\ref{thm:criterion}, $\D$ is claw-free.
		
		\emph{Necessity.} Let $K=G\times H$ and assume (A$'$), (B$'$), (C$'$) hold in $K$ (Theorem~\ref{thm:criterion}); coprimality of $|G|$ and $|H|$ is never used.
		
		\emph{Step 1.} Let $p\in\pi(G)$, $q\in\pi(H)$, $g\in G$ of order $p$, $h\in H$ of order $q$. Then $|(g,h)|=\mathrm{lcm}(p,q)$, which equals $pq$ if $p\ne q$. By (A$'$), $p=q$ or $\{p,q\}=\{2,3\}$.
		
		\emph{Step 2.} Suppose some prime $r\ge5$ lies in $\pi(G)$. By Step 1 every $q\in\pi(H)$ equals $r$; applying Step 1 again to each $p\in\pi(G)$ and $q=r$ gives $p=r$. So $\pi(G)=\pi(H)=\{r\}$, which is (1). The same holds if $r\in\pi(H)$. Otherwise $\pi(G),\pi(H)$ are non-empty subsets of $\{2,3\}$; if both equal $\{2\}$ or both equal $\{3\}$ we are in (1). Up to swapping $G$ and $H$ four cases remain.
		
		\emph{(a) $G$ a $2$-group, $H$ a $3$-group.} For $g\ne1$ of order $2^\alpha$ and $h\ne1$ of order $3^\beta$, $|(g,h)|=2^\alpha3^\beta$, so $\alpha=\beta=1$ by (A$'$): $G$ has exponent $2$ and $H$ exponent $3$. If $|G|\ge4$, take distinct involutions $g_1,g_2$ and $h$ of order $3$; then $(1,h)$ commutes with $(g_1,1)$ and $(g_2,1)$, contradicting (B$'$). So $|G|=2$. If $|H|\ge9$, the central involution $(u,1)$ centralizes the $|H|-1\ge8$ elements $(1,h)$ of order $3$, contradicting (C$'$). So $G=C_2$, $H=C_3$: (2).
		
		\emph{(b) $G$ a $2$-group, $\pi(H)=\{2,3\}$.} Let $a\in H$ have order $3$. For $g\in G$ of order $2^k$, $|(g,a)|=3\cdot2^k$, so $k\le1$: $G$ has exponent $2$. The element $(1,a)$ centralizes all $(g,1)$, $g\ne1$, so (B$'$) gives $|G|-1\le1$ and $G=\gen u\cong C_2$. The involutions of $C_K((1,a))=C_2\times C_H(a)$ include $(u,1)$ and $(1,v)$ for each involution $v\in C_H(a)$; by (B$'$), $C_H(a)$ has no involution, for every $a$ of order $3$. Since $(u,1)$ is central, (C$'$) shows $H$ has at most two elements of order $3$, hence exactly two, which form a unique subgroup $A$ of order $3$. $H$ has no element $y$ of order $9$ (else $|(u,y)|=18$). So a Sylow $3$-subgroup of $H$ has all its non-identity elements of order $3$ lying in $A$, hence equals $A$; thus $A\trianglelefteq H$ and $H=A\rtimes T$ with $T\ne1$ a Sylow $2$-subgroup. The kernel $C_T(a)$ of the action of $T$ on $A$ is a $2$-subgroup of $C_H(a)$, hence trivial. So $T$ embeds in $\Aut(A)\cong C_2$, $|T|=2$ and $T$ inverts $a$: $H\cong\Sym_3$ and $K\cong D_{12}$: (3).
		
		\emph{(c) $G$ a $3$-group, $\pi(H)=\{2,3\}$.} Let $g\in G$ have order $3$. Then $(g,1)$ centralizes every $(1,v)$, $v\in H$ an involution, so by (B$'$) $H$ has exactly one involution $v_0$, which is therefore central in $H$, and $(1,v_0)$ is central in $K$. By (C$'$) it centralizes at most two elements of order $3$; but $(g,1),(g^{-1},1),(1,a),(1,a^{-1})$ are four ($a\in H$ of order $3$). Contradiction.
		
		\emph{(d) $\pi(G)=\pi(H)=\{2,3\}$.} Let $a\in G$, $b\in H$ have order $3$. As in (c), (B$'$) at $(a,1)$ shows that $H$ has a unique involution $v_0$ and that $C_G(a)$ has no involution (the allowed one is already $(1,v_0)$). Symmetrically $G$ has a unique involution $u_0$, which is central in $G$ and so lies in $C_G(a)$: a contradiction.
		
		The last assertion: let $k\ge3$ and put $A=G_1\times\cdots\times G_{k-1}$, $B=G_k$. If $\D(A\times B)$ is claw-free and non-empty, Theorem~\ref{thm:direct} gives $\{A,B\}\in\{\{C_2,C_3\},\{C_2,\Sym_3\}\}$. But $A$ is a direct product of at least two non-trivial groups, whereas $C_2$, $C_3$ and $\Sym_3$ are directly indecomposable, a contradiction.
	\end{proof}
	
	\subsection{Semidirect products}
	
	\begin{proposition}\label{prop:semidirect-general}
		Let $X=N\rtimes K$ with $N,K\ne1$, let $K_0=C_K(N)$ be the kernel of the action, and assume $\D(X)$ is claw-free.
		\begin{enumerate}[label=\textup{(\alph*)}]
			\item $\D(N)$ and $\D(K)$ are claw-free.
			\item If $K_0\ne1$, then $N$ and $K_0$ are $p$-groups for the same prime $p$, or $\{N,K_0\}=\{C_2,C_3\}$ or $\{C_2,\Sym_3\}$.
			\item If $\D(X)\ne\emptyset$ and $N$ is a $2$-group, then $N=C_2$ and $X\cong C_2\times K$ with $K\in\{C_3,\Sym_3\}$. If $\D(X)\ne\emptyset$ and $N$ is a $3$-group, then $|N|=3$.
		\end{enumerate}
	\end{proposition}
	
	\begin{proof}
		(a) $N,K\le X$ and Lemma~\ref{lem:subgroups}. (b) $N\times K_0\le X$, so $\D(N\times K_0)$ is claw-free; apply Theorem~\ref{thm:direct}. (c) $N\le O_2(X)$ (respectively $O_3(X)$), so Lemma~\ref{lem:O2O3} gives $|N|\le2$ (respectively $\le3$). A normal subgroup of order $2$ is central, so $X=C_2\times K$ and Theorem~\ref{thm:direct}, with $\D(X)\ne\emptyset$, gives $K\in\{C_3,\Sym_3\}$ (note that $K$ is non-trivial and not a $2$-group, since otherwise $X$ is a $2$-group).
	\end{proof}
	
	\begin{theorem}\label{thm:semidirect-soluble}
		Let $X=N\rtimes K$ be soluble and non-EPPO, with $N,K\ne1$. Then $\D(X)$ is claw-free if and only if $X$ is one of the groups in Theorem~\textup{\ref{thm:soluble}}; for these groups the decompositions $X=N\rtimes K$ with $N\trianglelefteq X$ and $K$ a complement, $N,K\ne1$, are exactly the following.
		\begin{enumerate}[label=\textup{(\arabic*)}]
			\item $X=C_6$: $(C_2,C_3)$ or $(C_3,C_2)$ (trivial action).
			\item $X=D_{12}$: $(C_3,C_2^2)$, $(C_6,C_2)$, $(\Sym_3,C_2)$, $(C_2,\Sym_3)$.
			\item $X=\Dic_3$: only $(C_3,C_4)$.
			\item $X=M\rtimes K'$ Frobenius with $M$ an abelian $r$-group ($r\ge5$) and $K'\in\{C_6,\Dic_3\}$:
			\begin{enumerate}[label=\textup{(4\alph*)}]
				\item $N=M$, $K=K'$;
				\item $N=M_0$, $K=M_1\rtimes K'$, where $M=M_0\oplus M_1$ with $M_0,M_1\ne1$ both $K'$-invariant;
				\item $N=M\rtimes L$ with $L\trianglelefteq K'$ and $K$ a complement of $L$ in $K'$: $(M\rtimes C_2,C_3)$ and $(M\rtimes C_3,C_2)$ if $K'=C_6$, and $(M\rtimes C_3,C_4)$ if $K'=\Dic_3$.
			\end{enumerate}
		\end{enumerate}
	\end{theorem}
	
	\begin{proof}
		By Theorem~\ref{thm:soluble}, $X$ is $C_6$, $D_{12}$, $\Dic_3$ or a Frobenius group as in (4); it remains to list the pairs $(N,K)$ with $N\trianglelefteq X$ non-trivial proper and $K$ a complement.
		
		$C_6$ has normal subgroups $C_2,C_3$, each with a (unique) complement. The non-trivial proper normal subgroups of $D_{12}$ are $C_2=Z(D_{12})$, $C_3$, $C_6$ and the two subgroups $\cong\Sym_3$; their complements are $\Sym_3$, $C_2^2$, $C_2$, $C_2$ respectively. In $\Dic_3$ the non-trivial proper normal subgroups are $C_2$, $C_3$, $C_6$ (the Sylow $2$-subgroups $C_4$ are not normal, and the only subgroup of index $2$ is $C_6$). $\Dic_3$ has a unique involution $z$, which lies in every subgroup of even order. A complement of $C_2$ or of $C_6$ would have order $6$ or $2$, hence contain $z$, so would meet the normal subgroup non-trivially. Only $C_3$ has a complement, namely $C_4$.
		
		Let $X=M\rtimes K'$ be Frobenius. We use the standard fact that a normal subgroup $N$ of a Frobenius group with kernel $M$ satisfies $N\le M$ or $M\le N$ (see Isaacs, \emph{Finite Group Theory}, Ch.~6), and that, $X$ being soluble, all $r'$-Hall subgroups are conjugate (Hall; for the subgroup $K''$ below, Schur--Zassenhaus applied to $M\trianglelefteq MK''$).
		
		\emph{$1\ne N\le M$.} Since $NK=X$ and $N\le M$, $KM=X$. Put $M_1=K\cap M$, the (normal) Sylow $r$-subgroup of $K$; then $|M_1|=|K|/|K'|=|M:N|$, and $N\cap M_1=1$ gives $M=N\oplus M_1$. Write $K=M_1\rtimes K''$ with $K''$ an $r'$-Hall subgroup of $X$; by the Schur--Zassenhaus theorem $K''$ is conjugate to $K'$, and replacing $K$ by a conjugate (which is again a complement of the normal subgroup $N$) we may assume $K''=K'$. Then $M_1$ is normalized by $K'$, i.e.\ $K'$-invariant, and $N$ is $K'$-invariant because it is normal in $X$. Conversely, if $M=M_0\oplus M_1$ with $M_0,M_1$ $K'$-invariant, then $M_0\trianglelefteq X$ and $M_1K'$ is a complement. This gives (4a) when $M_1=1$ and (4b) otherwise.
		
		\emph{$M\le N$.} Then $N=M\rtimes L$, where $L=N\cap K'$ is normal in $K'$ (as $N/M\trianglelefteq X/M\cong K'$). Since $K\cap N=1$ and $NK=X$, $|K|=|K'|/|L|$ is coprime to $r$, so $K$ is conjugate into $K'$ and we may assume $K\le K'$, $K'=LK$, $L\cap K=1$. For $K'=C_6$ the normal subgroups $L$ of $K'$ with a complement are $1,C_2,C_3$ (and $C_6$, giving $K=1$, which is excluded), giving (4a) and (4c). For $K'=\Dic_3$: $L=C_2$ has no complement (a subgroup of order $6$ is a $C_6$ and contains the unique involution), $L=C_6$ has none (the complement would be generated by the unique involution, which lies in $C_6$), and $L=C_3$ has complement $C_4$.
		
		Conversely all listed pairs are decompositions. 
	\end{proof}
	
	\begin{remark}\label{rem:semidirect}
		(i) The group $\Dic_3$ is not a direct product of two non-trivial groups: its only decomposition is $C_3\rtimes C_4$ with non-trivial action. In case (4b), $K'$ acts fixed-point-freely on $M$, hence on each summand $M_0,M_1$; so if $K'=\Dic_3$ neither summand is cyclic (the automorphism group of a cyclic group is abelian), and (4b) needs $M$ to have at least four generators.
		(ii) If $X$ is EPPO then $\D(X)=\emptyset$, so every EPPO semidirect product is trivially claw-free; classifying them is the classification of EPPO groups.
		(iii) Non-soluble examples: $\Sym_5=\Alt_5\rtimes C_2$ and $\PGL_2(7)=\PSL_2(7)\rtimes C_2$ (the involutions of $\PGL_2(7)$ with determinant a non-square lie outside $\PSL_2(7)$, so the extension splits). A classification of the non-soluble claw-free semidirect products is open; Proposition~\ref{prop:semidirect-general} gives necessary conditions only.
		(iv) The decompositions in Theorem~\ref{thm:semidirect-soluble} are checked in GAP for all non-EPPO claw-free groups of order $\le200$, orders $128$ and $192$ excluded ($C_6$, $D_{12}$, $\Dic_3$, $C_7\rtimes C_6$, $C_{13}\rtimes C_6$, $C_{19}\rtimes C_6$, $C_{31}\rtimes C_6$, $C_5^2\rtimes C_6$), and, with orders $256$ and $384$ excluded, for orders $201$--$400$ (the additional groups being $C_{37}\rtimes C_6$, $C_{43}\rtimes C_6$, $C_{49}\rtimes C_6$, $C_7^2\rtimes C_6$, $C_{61}\rtimes C_6$, $C_5^2\rtimes\Dic_3$, and the non-soluble $\PGL_2(7)$); the direct product theorem for all products with $|G|,|H|\ge2$, $|G||H|\le400$, $|G|,|H|\le48$.
	\end{remark}

	\section*{Acknowledgments}
	The author gratefully acknowledges the financial support provided by the Anusandhan National Research Foundation (ANRF), India, under the National Postdoctoral Fellowship Scheme (File No.\ PDF/2025/006751). The author also acknowledges the excellent research facilities provided by the Department of Mathematics, IIT Hyderabad.

	\section*{Statements and Declarations}
	
	\noindent\textbf{Competing Interests.} The authors have no competing interests to declare that are
	relevant to the content of this article.
	
	\section*{Data Availability}
	
	Data sharing is not applicable to this article, as no datasets were used or generated during this
	study.

\end{document}